\documentclass{llncs}
\usepackage{epsfig}

\usepackage{makeidx}  % allows for indexgeneration
\usepackage{amsmath} \usepackage{amssymb, latexsym}

\begin{document}

\frontmatter

\def\supp{{\mbox{supp}}}

\title{Cross-Intersecting Families of Vectors}
\titlerunning{Cross-Intersecting Families of Vectors}

\author{J\'anos Pach\inst{1}
\and G\'abor  Tardos\inst{2}}
\authorrunning{J. Pach and G. Tardos}

\institute{\'{E}cole Polytechnique F\'ed\'erale de Lausanne\\
\email{pach@renyi.hu}
\and
R\'enyi Institute of Mathematics,
Budapest\\
and Center for Discrete Mathematics, Zhejiang Normal University\\
\email{tardos@renyi.hu}}
\maketitle

\begin{abstract}
Given a sequence of positive integers $p=(p_1,\dots,p_n)$, let $S_p$ denote
the family of all sequences of positive integers $x=(x_1,\ldots,x_n)$ such that
$x_i\leq p_i$ for all $i$. Two families of sequences (or vectors),
$A,B\subseteq S_p$, are said to be {\em $r$-cross-intersecting} if no matter
how we select $x\in A$ and $y\in B$, there are at least $r$ distinct indices
$i$ such that $x_i=y_i$. We determine the maximum
value of $|A|\cdot|B|$ over all pairs of $r$-cross-intersecting families and characterize
the extremal pairs for $r\ge 1$, provided that $\min p_i>r+1$. The case $\min p_i\le r+1$
is quite different. For this case, we have a conjecture, which we can verify
under additional assumptions. Our results generalize and strengthen several
previous results by Berge, Borg, Frankl, F\"uredi, Livingston, Moon, and Tokushige, and answers
a question of Zhang.
\end{abstract}

\let\thefootnote\relax\footnotetext
{J. Pach is supported by OTKA under ERC
    projects GraDR and ComPoSe 10-EuroGIGA-OP-003, and by Swiss National
    Science Foundation Grants 200020-144531 and 200021-137574.
    \\
    G. Tardos is
    supported by
    OTKA grant NN-102029, the ``Lend\"ulet'' project
    of the Hungarian Academy of Sciences and by EPFL.}

\section{Introduction}

The Erd\H os-Ko-Rado theorem~\cite{EKR61} states that for $n\geq 2k$, every family  of pairwise intersecting $k$-element subsets of an $n$-element set consists of at most ${n-1\choose k-1}$ subsets, as many as the star-like family of all subsets containing a fixed element of the underlying set. This was the starting point of a whole new area within combinatorics: extremal set theory; see~\cite{GK78}, \cite{Bol86}, \cite{DeF83}, \cite{F95}. The Erd\H os-Ko-Rado theorem has been extended and generalized to other structures: to multisets, divisors of an integer, subspaces of a vector space, families of permutations, etc. It was also generalized to ``cross-intersecting" families, i.e., to families $A$ and $B$ with the property that every element of $A$ intersects all elements of $B$; see Hilton~\cite{Hi77}, Moon~\cite{Mo82}, and Pyber~\cite{Py86}.

\smallskip

For any positive integer $k$, we write $[k]$ for the set $\{1,\dots,k\}$. Given a sequence of positive integers $p=(p_1,\dots,p_n)$, let
$$S_p=[p_1]\times\cdots\times[p_n]=\{(x_1,\dots,x_n)\; :\; x_i\in[p_i]\hbox{ for }i\in[n]\}.$$
We will refer to the elements of $S_p$ as {\em vectors}. The {\em Hamming
  distance} between the vectors $x,y\in S_p$ is $|\{i\in[n]\;:\; x_i\ne y_i\}|$
and is denoted by $d(x,y)$. Let $r\ge1$ be an
integer. Two vectors $x,y\in
S_p$ are said to be {\em $r$-intersecting} if $d(x,y)\le n-r$. (This term
originates in the
observation that if we represent a vector $x=(x_1,\dots,x_n)\in S_p$ by the
set $\{(i,x_i)\; :\; i\in[n]\}$, then $x$ and $y\in S_p$ are $r$-intersecting if
and only if the sets representing them have at least $r$ common elements.)
Two families $A,B\subseteq S_p$ are {\em $r$-cross-intersecting}, if every
pair $x\in A$, $y\in B$ is $r$-intersecting. If $(A,A)$ is an
$r$-cross-intersecting pair, we say $A$ is {\em $r$-intersecting}. We simply say {\em
  intersecting} or {\em cross-intersecting} to mean $1$-intersecting or
$1$-cross-intersecting, respectively.
\smallskip

The investigation of the maximum value for $|A|\cdot|B|$ for cross-intersecting pairs of families $A,B\subseteq S_p$ was initiated by Moon~\cite{Mo82}. She proved, using a clever induction argument, that in the special case when  $p_1=p_2=\dots=p_n=k$ for some $k\ge 3$, every cross-intersecting pair $A,B\subseteq S_p$ satisfies
$$|A|\cdot|B|\le k^{2n-2},$$
with equality if and only if $A=B=\{x\in S_p\; :\; x_i=j\}$, for some $i\in[n]$ and $j\in[k]$.  In the case $A=B$, Moon's theorem had been discovered by Berge~\cite{Be74}, Livingston~\cite{Liv79}, and Borg~\cite{Bo08}. See also Stanton~\cite{St80}. In his report on Livingston's paper, published in the {\em Mathematical Reviews}, Kleitman gave an extremely short proof for the case $A=B$, based on a shifting argument. Zhang~\cite{Zh13} established a somewhat weaker result, using a generalization of Katona's circle method~\cite{K72}. Note that for $k=2$, we can take $A=B$ to be any family of $2^{n-1}$ vectors without containing a pair $(x_1,\ldots,x_n), (y_1,\ldots,y_n)$ with $x_i+y_i=3$ for every $i$. Then $A$ is an intersecting family with $|A|^2=2^{2n-2}$, which is not of the type described in Moon's theorem.

Moon also considered $r$-cross-intersecting pairs in $S_p$ with $p_1=p_2=\dots=p_n=k$ for some $k > r+1$, and characterized all pairs for which $|A|\cdot|B|$ attains its maximum, that is, we have
$$|A|\cdot|B|= k^{2(n-r)}.$$
The assumption $k>r+1$ is necessary. See Tokushige~\cite{To13}, for a somewhat weaker result, using algebraic techniques.
%\gabor{HIANYZIK MEG HIVATKOZAS???}

\smallskip

Zhang~\cite{Zh13} suggested that Moon's results may be extended to arbitrary
sequences of positive integers $p=(p_1,\dots,p_n)$. The aim of this note is
twofold: (1) to establish such an extension under the assumption $\min_i p_i > r+1$,
and (2) to formulate a conjecture that covers essentially all
other interesting cases. We verify this conjecture in several special cases.

We start with the special case $r=1$, which has also been settled independently by Borg~\cite{Bo14}, using different techniques.

\medskip
\noindent
{\bf Theorem 1.} {\em Let $p=(p_1,\dots,p_n)$ be a sequence positive integers
  and let $A,B\subseteq S_p$ form a pair of cross-intersecting families of
  vectors.

  We have $|A|\cdot|B|\le|S_p|^2/k^2$, where $k=\min_ip_i$. Equality holds for
  the case $A=B=\{x\in S_p\; :\; x_i=j\}$, whenever $i\in[n]$ satisfies $p_i=k$
  and $j\in[k]$. For $k\ne2$, there are no other extremal cross-intersecting
  families.}
\medskip

We say that a coordinate $i\in[n]$ is {\em irrelevant} for a set $A\subseteq
S_p$ if, whenever two elements of $S_p$ differ only in coordinate $i$ and $A$
contains one of them, it also contains the other. Otherwise, we say that $i$
is {\em relevant} for $A$.

Note that no coordinate $i$ with $p_i=1$ can be relevant for any family. Each such
coordinate forces an intersection between every pair of vectors. So, if we
delete it, every $r$-cross-intersecting pair becomes $(r-1)$-cross-intersecting.
Therefore, from now on we will always assume that we have
$p_i\ge2$ for every $i$.

We call a sequence of integers
$p=(p_1,\dots,p_n)$ a {\em size vector} if $p_i\ge2$ for all $i$. The {\em
length} of $p$ is $n$.
We say that an $r$-cross-intersecting pair $A,B\subseteq
S_p$ is {\em maximal} if it maximizes the value $|A|\cdot|B|$.

\smallskip
Using this notation and terminology, Theorem~1 can be rephrased as follows.

\medskip
\noindent{\bf Theorem 1'.}
{\em Let $p=(p_1,\dots,p_n)$ be a sequence of integers with $k=\min_ip_i>2$.

For any maximal pair of cross-intersecting families,
$A,B\subseteq S_p$, we have $A=B$, and there is a single coordinate which is
relevant for $A$. The relevant coordinate $i$ must satisfy $p_i=k$.}
\medskip

See Section~\ref{cr} for a complete characterization of maximal
cross-intersecting pairs in the $k=2$ case. Here we mention that only the
coordinates with $p_i=2$ can be relevant for them, but
%\janos{We have to argue/}
for certain pairs, {\em all} such coordinates are relevant
simultaneously. For example, let $n$ be odd, $p=(2,\ldots,2)$, and let $A=B$
consist of all vectors in $S_p$ which have at most $\lfloor n/2\rfloor$
coordinates that are $1$. This makes $(A,B)$ a maximal cross-intersecting pair.

\smallskip

Let $T\subseteq[n]$ be a subset of the coordinates, let $x_0\in S_p$ be an
arbitrary vector, and let $k$ be an integer satisfying $0\le k\le|T|$. The {\em Hamming ball} of
radius $k$ around $x_0$ in the coordinates $T$ is defined as the family
$$B_k=\{x\in S_p\; :\;|\{i\in T\; :\; x_i\ne(x_0)_i\}|\le k\}.$$
Note that the pair $(B_k,B_l)$ is
$(|T|-k-l)$-cross-intersecting. We use the word {\em ball} to refer to any
Hamming ball without specifying its center, radius
or its set of coordinates. A Hamming ball of radius $0$ in coordinates $T$ is
said to be obtained by {\em fixing} the coordinates in $T$.

For the proof of Theorem~1, we need the following statement, which will be established by induction on $n$, using the idea in~\cite{Mo82}.

\medskip
\noindent
{\bf Lemma 2.} {\em Let $1\le r<n$, let $p=(p_1,\dots,p_n)$ be a size
  vector satisfying $3\le p_1\le p_2\le\cdots\le p_n$ and let
  $A,B\subseteq S_p$ form a pair of $r$-cross-intersecting families. If
  $$\frac{2}{p_{r+1}}+\sum_{i=1}^r\frac{1}{p_i}\le1,$$
  then
  $|A|\cdot|B|\le\prod_{i=r+1}^np_i^2$. In case of equality, we have
  $A=B$ and this family can be obtained by fixing $r$ coordinates in
  $S_p$.}
\medskip

By fixing any $r$ coordinates, we obtain a ``trivial'' $r$-intersecting
family $A=B\subseteq S_p$. As was observed by Frankl and F\"uredi \cite{FF80}, not all maximal size $r$-intersecting families can be obtained in this way, for certain size vectors. They considered size vectors $p=(k,\ldots,k)$
with $n\ge r+2$ coordinates, and noticed that a Hamming ball of radius $1$ in $r+2$ coordinates is $r$-intersecting. Moreover, for $k\le r$, this family is strictly larger than the trivial $r$-intersecting family. See also~\cite{AhK98}.

On the other hand, as was mentioned before, for $k\ge r+2$, Moon \cite{Mo82} proved that among all $r$-intersecting families, the trivial ones are maximal.

\smallskip

This leaves open only the case $k=r+1$, where the trivial $r$-intersecting
families and the radius $1$ balls in $r+2$ coordinates have precisely the same
size. We believe that in this case there are no larger $r$-intersecting
families. For $r=1$, it can be and has been easily verified (and follows, for
example, from our Theorem~1, which deals with the asymmetric case, when $A$
and $B$ do not necessarily coincide). Our Theorem~7 settles the problem also
for $r>3$. The intermediate cases $r=2$ or $3$ are still open, but they could possibly be handled by computer search.

\smallskip

Therefore, to characterize maximal size $r$-intersecting families $A$ or
maximal $r$-cross-intersecting pairs of families $(A,B)$ for {\em all} size
vectors, we cannot restrict ourselves to fixing $r$ coordinates. We make the
following conjecture that can roughly be considered as a generalization of the
Frankl-F\"uredi conjecture \cite{FF80} that has been proved by Frankl and Tokushige \cite{FT99}. The
generalization is twofold: we consider $r$-cross-intersecting pairs rather
than $r$-intersecting families and we allow arbitrary size vectors not just
vectors with all-equal coordinates.
\medskip

\noindent
{\bf Conjecture 3.} {\em If $1\le r\le n$ and $p$ is a size
  vector of length $n$, then there exists a maximal pair of
  $r$-cross-intersecting families $A,B\subseteq S_p$, where $A$ and $B$ are
  balls. If we further have $p_i\ge3$ for all $i\in[n]$, then all maximal
  pairs of $r$-cross-intersecting families consist of balls.}
\medskip

Note that the $r=1$ special case of Conjecture~3 is established by
Theorem~1. Some further special cases of the conjecture are settled in Theorem~7.

It is not hard to narrow down the range of possibilities for maximal
$r$-cross-intersecting
pairs that are formed by two balls, $A$ and $B$. In fact, the following theorem
implies that all such pairs are determined up to isomorphism by the
number of relevant coordinates. Assuming that Conjecture~3 is true, finding
$\max|A|\cdot|B|$ for $r$-cross-intersecting families $A,B\subseteq S_p$ boils
down to making numeric comparisons for pairs of balls obtained by various radii. In case $p_i\ge3$ for all $i$ (and still assuming
Conjecture~3), the same process also finds all maximal $r$-cross-intersecting
pairs.
\medskip

\noindent
{\bf Theorem 4.} {\em Let $1\le r\le n$ and let $p=(p_1,\dots,p_n)$ be a size
  vector. If $A,B\subseteq S_p$ form a maximal pair of $r$-cross-intersecting
  families, then either of them determines the other. In particular, $A$ and
  $B$ have the same set of relevant coordinates. Moreover, if $A$ is a ball of radius
  $l$ around $x_0\in S_p$ in a set of coordinates $T\subseteq[n]$, then
  $|T|\ge l+r$, and $B$ is a ball of radius $|T|-l-r$ around $x_0$ in the same set
  of the coordinates. Furthermore, we have $p_i\le p_j$ for every $i\in T$ and
  $j\in[n]\setminus T$, and the radii of the balls differ by at most 1, that is,
  $\big||T|-2l-r\big|\le1$.
}
\medskip

Note that if $A=B$ for a maximal pair $(A,B)$ of $r$-cross-intersecting
families, then $A$ is also a maximal size $r$-intersecting family. This is the
case, in particular, if $A$ and $B$ are balls of equal radii. However, for many size
vectors, no maximal $r$-cross-intersecting pair consists of $r$-intersecting families, as the maximal $r$-cross-intersecting
pairs are often formed by balls whose radii differ by one. For example, for the size
vector $p=(3,3,3,3,3)$, the largest $4$-intersecting family $C$ is obtained by
fixing four coordinates, while the maximal $4$-cross-intersecting pair is formed
by a singleton $A=\{x\}$ and a ball $B$ of radius $1$ around $x$ in all 
coordinates. Here we have $|A|\cdot|B|=11>|C|^2=9$.

As we have indicated above, we have been unable to prove Conjecture~3 in its full generality, but we were able to verify it
in several interesting special cases. We will proceed in two steps. First we argue, using {\em entropies}, that the number of relevant coordinates in a maximal $r$-cross-intersecting family is bounded. Then we apply combinatorial methods to prove the conjecture under the assumption that the number of relevant coordinates is small.

In the case where there are many relevant coordinates for a pair of maximal
$r$-cross-intersecting families, we use entropies to bound the size of the
families and to prove

\medskip
\noindent
{\bf Theorem 5.} {\em Let $1\le r\le n$, let $p=(p_1,\dots,p_n)$ be a size
  vector, let $A,B\subseteq S_p$ form a maximal pair of $r$-cross-intersecting
  families, and let $T$ be the set of coordinates that are relevant for $A$ or
  $B$. Then neither the size of $A$ nor the size of $B$ can exceed
$$\frac{|S_p|}{\prod_{i\in T}(p_i-1)^{1-2/p_i}}.$$
}
\medskip

We use this theorem to bound the number of relevant coordinates $i$ with
$p_i>2$. The number of relevant coordinates $i$ with $p_i=2$ can be unbounded; see Section~5.

\medskip
\noindent
{\bf Theorem 6.} {\em Let $1\le r\le n$, let $p=(p_1,\dots,p_n)$ be a size
  vector, and let $A,B\subseteq S_p$ form a maximal pair of
  $r$-cross-intersecting families.

  For the set of coordinates $T$ relevant for
  $A$ or $B$, we have $$\prod_{i=1}^rp_i\ge\prod_{i\in T}(p_i-1)^{1-2/p_i},$$
  which implies that $|\{i\in T\; :\; p_i>2\}|<5r$.}
\medskip

We can characterize the maximal $r$-cross-intersecting pairs for
all size vectors $p$ satisfying $\min p_i>r+1$, and in many other cases.
\medskip

\noindent
{\bf Theorem 7.} {\em Let $1\le r\le n$, let $p=(p_1,\dots,p_n)$ be a size
  vector with $p_1\le p_2\le\cdots\le p_n$, and let $A,B\subseteq
  S_p$ form a pair of $r$-cross-intersecting families.
\begin{enumerate}
\item  If $p_1>r+1$,
  we have $|A|\cdot|B|\le\prod_{i=r+1}^np_i^2$. In case of equality,
  $A=B$ holds and this family can be obtained by fixing $r$ coordinates in $S_p$.
\item If $p_1=r+1>4$, we have $|A|\cdot|B|\le\prod_{i=r+1}^np_i^2$. In case of equality,
  $A=B$ holds and this family can be obtained either by fixing
  $r$ coordinates in $S_p$ or by taking a Hamming ball of radius $1$ in
  $r+2$ coordinates $i$, all satisfying $p_i=r+1$.
\item There is a function $t(r)=r/2+o(r)$ such that if $p_1\ge t(r)$ and
  $(A,B)$ is a maximal $r$-cross-intersecting pair,
  then the families $A$ and $B$ are balls in at most $r+3$
  coordinates.
\end{enumerate}}
\medskip

The proof of Theorem~7 relies on the following result.

\bigskip
\bigskip
\noindent
{\bf Theorem 8.} {\em Let $1\le r\le n$ and let
  $p$ be a size vector of length $n$.
\begin{enumerate}
\item If there exists a maximal pair of $r$-cross-intersecting families in
  $S_P$ with at most $r+2$ relevant coordinates,
  then there exists such a pair consisting of balls.
\item If $p_i>2$ for all $i\in[n]$ and
  $A,B\subseteq S_p$ form a maximal pair of $r$-cross-intersecting families
  with at most $r+3$ relevant coordinates, then $A$ and $B$ are balls.
\end{enumerate}}
\medskip

With an involved case analysis, it may be possible to extend Theorem~8 to
pairs with more relevant coordinates. Any such an improvement would carry
over to Theorem~7.

All of our results remain meaningful in the symmetric case where $A=B$. For instance, in
this case, Theorem~1 (also proved by Borg \cite{Bo14}) states that every intersecting family $A\subseteq S_p$
has at most $|S_p|/k$ members, where $k=\min_ip_i$. In case $k>2$,
equality can be achieved only by fixing some coordinate $i$ with $p_i=k$. Note
that in the case $A=B$ (i.e., $r$-intersecting families) the exact maximum
size is known for size vectors $(q,\ldots,q)$, \cite{FT99}.

\section{Proofs of Theorems 4 and 1}

First, we verify Theorem~4 and a technical lemma (see Lemma~9 below) which generalizes the corresponding result in~\cite{Mo82}. Our proof is slightly simpler. Lemma~9 will enable us to deduce Lemma~2, the main ingredient of the proof of Theorem~1, presented at the end of the section.

\medskip
\noindent
{\em Proof of Theorem 4.} The first statement is self-evident: if $A,B\subseteq
S_p$ form a maximal pair of $r$-cross-intersecting families, then
$$B=\{x\in S_p\; :\; x\hbox{ $r$-intersects $y$ for all }y\in A\}.$$
If a coordinate is irrelevant for $A$, then it is also irrelevant for $B$ defined by this formula. Therefore, by symmetry, $A$ and $B$ have the same set of relevant coordinates.

If $A$ is the Hamming ball around $x_0$ of radius $l$ in
coordinates $T$, then we have $B=\emptyset$ if $|T|<l+r$, which is not
possible for a maximal cross-intersecting family. If $|T|\ge l+r$, we obtain
the ball claimed in the theorem. For every $i\in T$, $j\in[n]\setminus T$,
consider the set $T'=(T\setminus\{i\})\cup\{j\}$ and the Hamming balls $A'$
and $B'$ of radii $l$ and $|T|-l-r$ around $x_0$ in the coordinates
$T'$. These balls form an $r$-cross-intersecting pair and in case $p_i>p_j$, we
have $|A'|>|A|$ and $|B'|>|B|$, contradicting the maximality of the pair
$(A,B)$.

Finally let $B_l$ be a ball of radius $l$ around some fixed vector $x$ in
a fixed set $T$ of coordinates. We claim that the size $|B_l|$ of these balls
is strictly log-concave, that is, we have
$$|B_l|^2>|B_{l-1}|\cdot|B_{l+1}|$$
for $1\le l<|T|$. As balls around different centers have the same size, we can
represent the left-hand side as $|B_l|^2=|C|$, where
$$C=\{(y,z)\mid y,z\in S_p, d(x,y)\le l,d(y,z)\le l\}.$$
Similarly, the right-hand side can be
represented as $|B_{l-1}|\cdot|B_{l+1}|=|D|$ with
$$D=\{(y,z)\mid y,z\in S_p,d(x,y)\le l-1,d(y,z)\le l+1\}.$$

We say that two pairs $(y,z)$ and $(y',z')$ (all four terms from $S_p$)
are {\em equivalent} if $z=z'$ and for every $i\in[n]$ we have either
$y_i=y'_i$ or $y_i,y'_i\in\{x_i,z_i\}$. Let us fix an equivalence class
$O$. For all $(y,z)\in O$, the element $z$ and some
coordinates $y_i$ of $y$ are {\em fixed}. We call the remaining coordinates $i$ {\em
open}. For an open coordinate $i$, the value of $y_i$ must be one of the
non-equal values $x_i$ or $z_i$. If $m$ denotes the number of open
coordinates in $O$, we have $|O|=2^m$. For a pair $(y,z)\in O$, the
distance $d(x,y)=d_1+d_2(y)$, where $d_1$ is the number of fixed
coordinates $i$ with $x_i\ne y_i$, while $d_2(y)$ is the number of open
coordinates $i$ with $x_i\ne y_i$. Note that $d_1$ is constant for all
elements of $O$, while $d_2(y)$ takes any value $j$ for $m\choose j$ members
of $O$. Similarly, we can write $d(y,z)=d_1+(m-d_2(y))$, as $y_i\ne z_i$ holds
for the same fixed coordinates where $x_i\ne y_i$, and $y_i$ is equal to
exactly one of $x_i$ and $z_i$ for an open coordinate $i$. Summarizing, we
have
$$C\cap O=\{(y,z)\in O\mid d_1+m-l\le d_2(y)\le l-d_1\},$$
$$D\cap O=\{(y,z)\in O\mid d_1+m-l-1\le d_2(y)\le l-d_1-1\}.$$
We claim that $|C\cap O|\ge|D\cap O|$. Indeed, if $l-d_1<m/2$, then $C\cap
O=D\cap O=\emptyset$. Otherwise, we have $|C\cap O|-|D\cap O|={m\choose
  l-d_1}-{m\choose l-d_1+1}\ge0$. Note also that equality holds only if
$l-d_1<m/2$ or $l-d_1>m$, in which cases $C$ and $D$ are disjoint from $O$ or
contain $O$, respectively.

As $C$ contains at least as many pairs from every equivalence class as $D$
does, we have $|C|\ge|D|$. Equality cannot hold for all equivalence classes, so
we have $|C|>|D|$, as claimed.

To finish the proof of the theorem, we need to verify that the pair
$(B_{l_1},B_{l_2})$ is not maximal $r$-cross-intersecting if $r=|T|-l_1-l_2$
and $|l_1-l_2|\ge2$. This follows from the log-concavity, because in case $l_1\ge
l_2+2$ the pair $(B_{l_1-1},B_{l_2+1})$ is also $r$-cross-intersecting and $|B_{l_1-1}|\cdot|B_{l_2+1}|>|B_{l_1}|\cdot|B_{l_2}|$.
\qed

\medskip

The following lemma will also be used in the proof of Theorem~5, presented in the next section.

\medskip

\noindent
{\bf Lemma 9.} {\em Let $1\le r\le n$, let $p=(p_1,\dots,p_n)$ be a size vector, and let $A,B\subseteq S_p$ form a maximal pair of $r$-cross-intersecting families.

If $i\in[n]$ is a relevant coordinate for $A$ or $B$, then there exists a value $l\in[p_i]$ such that
$$|\{x\in A\; :\; x_i\ne l\}|\le|A|/p_i,$$
$$|\{y\in B\; :\; y_i\ne l\}|\le|B|/p_i.$$
%Furthermore, for any $l'\in[p_i]\setminus\{l\}$ we have
%$$|\{x\in A\; :\; x_i=l'\}|\le|\{x\in A\;:\; x_i=l\}|/(p_i-1)^2,$$
%$$|\{x\in B\; :\; x_i=l'\}|\le|\{x\in B\; :\; x_i=l\}|/(p_i-1)^2.$$
}
\medskip

\noindent
{\em Proof.} Let us fix $r, n, p, i, A,$ and $B$ as in the
lemma. By Theorem~4, if a coordinate is irrelevant for $A$, then
it is also irrelevant for $B$ and vice versa.

In the case $n=r$, we have $A=B$ and this family must be a singleton, so the lemma is trivially true. From now on, we assume that $n>r$ and hence the notion of $r$-cross-intersecting families is meaningful for $n-1$ coordinates.

Let $q=(p_1,\ldots,p_{i-1},p_{i+1},\ldots,p_n)$. For any $l\in[p_i]$, let
$$A'_l=\{x\in A\; :\; x_i=l\},$$
$$B'_l=\{y\in B\; :\; y_i=l\},$$
and let $A_l$ and $B_l$ stand for the families obtained from $A'_l$ and $B'_l$, respectively, by dropping their $i$\/th coordinates. By definition, we have $A_l,B_l\subseteq S_q$, and $|A|=\sum_l|A_l|$ and $|B|=\sum_l|B_l|$. Furthermore, for any two distinct elements $l,m\in[p_i],$ the families $A_l$ and $B_m$ are $r$-cross-intersecting, since the vectors in $A'_l$ differ from the vectors in $B'_m$ in the $i$\/th coordinate, and therefore the $r$ indices where they agree must be elsewhere.

Let $Z$ denote the maximum product $|A^*|\cdot|B^*|$ of an
$r$-cross-intersecting pair $A^*,B^*\subseteq S_q$. We have
$|A_l|\cdot|B_m|\le Z$ for all
$l,m\in[p_i]$ with $l\ne m$. Adding an irrelevant $i$\/th coordinate to the
maximal $r$-cross-intersecting pair $A^*,B^*\subseteq S_q$, we obtain a pair
$A^{*\prime},B^{*\prime}\subseteq S_p$ with
$|A^{*\prime}|\cdot|B^{*\prime}|=p_i^2Z$. Using the maximality of $A$ and
$B$, we have $|A|\cdot|B|\ge p_i^2Z$. Let $l_0$ be chosen so as to maximize
$|A_{l_0}|\cdot |B_{l_0}|$, and let $c=|A_{l_0}|\cdot|B_{l_0}|/Z$.

Assume first that $c\le1$. Then we have
$$p_i^2Z\le|A|\cdot|B|=\sum_{l,m\in[p_i]}|A_l|\cdot|B_m|
\le\sum_{l,m\in[p_i]}Z=p_i^2Z.$$
Hence, we must have equality everywhere. This yields that $c=1$ and that $A_l$ and
$B_m$ form a maximal $r$-cross-intersecting pair for all
$l,m\in[p_i]$, $l\ne m$. This also implies that $|A_l|=|A_m|$ for $l,m\in[p_i]$, from where
the statement of the lemma follows, provided that $p_i=2$.

If $p_i\ge3$, then all families $A_l$ must be equal to one another, since one member in a maximal
$r$-cross-intersecting family determines the other, by Theorem~4. This contradicts our
assumption that the $i$\/th coordinate was relevant for $A$.

Thus, we may assume that $c>1$.
\smallskip

For $m\in[p_i]$, $m\ne l_0$, we have
$|A_{l_0}|\cdot|B_m|\le Z=|A_{l_0}|\cdot|B_{l_0}|/c$. Thus,
\begin{equation}\label{eq0}
|B_m|\le|B_{l_0}|/c,
\end{equation}
which yields that
$|B|=\sum_{m\in[p_i]}|B_m|\le(1+(p_i-1)/c)|B_{l_0}|$. By symmetry, we also
have
\begin{equation}\label{eq0b}
|A_m|\le|A_{l_0}|/c
\end{equation}
for $m\ne l_0$ and $|A|\le(1+(p_i-1)/c)|A_{l_0}|$. Combining these inequalities, we obtain
$$p_i^2Z\le|A|\cdot|B|\le(1+(p_n-1)/c)^2|A_{l_0}|\cdot|B_{l_0}|=(1+(p_i-1)/c)^2cZ.$$
We solve the resulting inequality $p_i^2\le c(1+(p_i-1)/c)^2$ for
$c>1$ and conclude that $c\ge(p_i-1)^2$. This inequality, together with Equations
(\ref{eq0}) and (\ref{eq0b}), completes the proof of Lemma~9. \qed

\medskip

\noindent
{\em Proof of Lemma 2.} We proceed by induction on $n$.

Let $A$ and $B$ form a maximal $r$-cross-intersecting pair. It is sufficient to
show that they have only $r$ relevant coordinates. Let us suppose that
the set $T$ of their relevant coordinates satisfies $|T|>r$, and choose a subset
$T'\subseteq T$ with $|T'|=r+1$. By Lemma~9, for every $i\in T'$ there exists
$l_i\in[p_i]$ such that the family $X_i=\{x\in B\; :\; x_i\ne l_i\}$ has cardinality $|X_i|\le|B|/p_i$.

If we assume that
$$\frac{2}{p_{r+1}}+\sum_{i=1}^r\frac{1}{p_i}<1$$
holds (with strict inequality), then this bound of $|X_i|$ would suffice.
In order to also be able to deal with the case
$$\frac{2}{p_{r+1}}+\sum_{i=1}^r\frac{1}{p_i}=1,$$
we show that $|X_i|=|B|/p_i$ is not possible. Considering the proof of
Lemma~9, equality here would mean that the families $A_l$ and $B_l$ (obtained
by dropping the $i$\/th coordinate from the vectors in the sets $\{x\in A\;
:\; x_i=l\}$ and $\{y\in B\; :\; y_i=l\}$, respectively) satisfy the following
condition: both $(A_{l_i},B_m)$ and
$(A_m,B_{l_i})$ should be maximal $r$-cross-intersecting pairs for all $m\ne
l_i$. By the induction hypothesis, this would imply that $A_{l_i}=B_m$ and
$A_m=B_{l_i}$, contradicting that $|A_m|<|A_{l_i}|$ and $|B_m|<|B_{l_i}|$ (see
(\ref{eq0}), in view of $c>1$). Therefore, we have $|X_i|<|B|/p_i$.

Let $C=\{x\in S_p\; :\; x_i=1\hbox{ for all }i\in[r]\}$ be the $r$-intersecting
family obtained by fixing $r$ coordinates in $S_p$.
In the family $D=B\setminus(\bigcup_{i\in T'}X_i)$, the coordinates in $T'$
are fixed. Thus, we have $$|D|\le\prod_{i\in[n]\setminus
T'}p_i\le\prod_{i=r+2}^np_i=|C|/p_{r+1}.$$
On the other hand, we have
$$|D|=|B|-\sum_{i\in T'}|X_i|>|B|(1-\sum_{i\in T'}1/p_i)\ge|B|(1-\sum_{i=1}^{r+1}1/p_i).$$
Comparing the last two inequalities, we obtain
$$|B|<\frac{|C|}{p_{r+1}(1-\sum_{i=1}^{r+1}1/p_i)}.$$
By our assumption on $p$, the denominator is at least $1$, so that we have
$|B|<|C|$. By symmetry, we also have $|A|<|C|$. Thus, $|A|\cdot|B|<|C|^2$
contradicting the maximality of the pair $(A,B)$. This completes the proof of
Lemma~2. \qed
\medskip

Now we can quickly finish the proof of Theorem~1.

\medskip
\noindent
{\em Proof of Theorem~1.} Notice that Lemma~2 implies Theorem~1, whenever
$k=\min_ip_i\ge3$. It remains to verify the statement for $k=1$ and $k=2$.
For $k=1$, it follows from the fact that all pairs of vectors in $S_p$ are
intersecting, thus the only maximal cross-intersecting pair is $A=B=S_p$.

Suppose next that $k=2$. For $x\in S_p$, let $x'\in S_p$ be defined
by $x'_i=(x_i+1\bmod p_i)$ for $i\in[n]$.
Note that $x\mapsto x'$ is a permutation of  $S_p$.
Clearly, $x$ and $x'$ are not intersecting, so we either
have $x\notin A$ or $x'\notin B$. As a
consequence, we obtain that $|A|+|B|\le|S_p|$, which, in turn, implies that
$|A|\cdot|B|\le|S_p|^2/4$, as claimed. It also follows that all maximal
pairs satisfy $|A|=|B|=|S_p|/2$. \qed

\section{Using entropy: Proofs of Theorems~5 and~6}

\noindent
{\em Proof of Theorem~5.} Let $r, n, p, A, B$ and $T$ be as in the theorem. Let
us write $y$ for a randomly and uniformly selected element of $B$.
Lemma~9 implies that, for each $i\in T$, there exists a value $l_i\in[p_i]$
such that
\begin{equation}\label{eq1}
Pr[y_i=l_i]\ge1-1/p_i.
\end{equation}

We bound the {\em entropy} $H(y_i)$ of $y_i$ from above by the entropy of the indicator variable of the event $y_i=l_i$ plus the contribution coming from the entropy of $y_i$ assuming $y_i\ne l_i$:
$$H(y_i)\le h(1-1/p_i)+(1/p_i)\log(p_i-1)=\log p_i-(1-2/p_i)\log(p_i-1),$$
where $h(z)=-z\log z-(1-z)\log(1-z)$ is the entropy function, and we used that $1-1/p_i\ge1/2$.

For any $i\in[n]\setminus T$, we use the trivial estimate $H(y_i)\le\log p_i$. By the subadditivity of the entropy, we obtain
$$\log|B|=H(y)\le\sum_{i\in[n]}H(y_i)\le\sum_{i\in T}(\log
p_i-(1-2/p_i)\log(p_i-1))+\sum_{i\in[n]\setminus T}\log p_i,$$
or, equivalently,
$$|B|\le\prod_{i\in T}\frac{p_i}{(p_i-1)^{1-2/p_i}}\prod_{i\in[n]\setminus
  T}p_i=\frac{|S_p|}{\prod_{i\in T}(p_i-1)^{1-2/p_i}}$$
as required. The bound on $|A|$ follows by symmetry and completes the proof of
the theorem. \qed
\medskip

Theorem~6 is a simple corollary of Theorem~5.

\medskip
\noindent
{\em Proof of Theorem~6.}
Fixing the first $r$ coordinates, we obtain the family
$$C=\{x\in S_p\; :\; x_i=1\hbox{ for all }i\in[r]\}.$$
This family is $r$-intersecting. Thus, by the maximality of the pair $(A,B)$, we have
\begin{equation}\label{eq3}
|A|\cdot|B|\ge|C|^2=\left(\prod_{i=r+1}^np_i\right)^2.
\end{equation}
Comparing this with our upper bounds on $|A|$ and $|B|$, we obtain the
inequality claimed in the theorem.

To prove the required bound on the number of relevant coordinates $i$ with
$p_i\ne2$, we assume that the coordinates are ordered, that is, $p_1\le
p_2\le\cdots\le p_n$. Applying the above estimate on $\prod_{i\in[r]}p_i$
and using $(p_i-1)^{1-2/p_i}>p_i^{1/5}$ whenever $p_i\ge3$, the theorem
follows. \qed

\section{Monotone families: Proofs of Theorems~8 and~7}

Given a vector $x\in S_p$, the set $\supp(x)=\{i\in[n] : x_i>1\}$ is called the
{\em support} of $x$. A family $A\subseteq S_p$ is said to be {\em monotone}, if for any
$x\in A$ and $y\in S_p$ satisfying $\supp(y)\subseteq\supp(x)$, we have $y\in A$.

For a family $A\subseteq S_p$, let us define its {\em support} as
$\supp(A)=\{\supp(x)\; :\; x\in A\}$. For a monotone family $A$, its support is
clearly subset-closed and it uniquely determines $A$, as $A=\{x\in
S_p\, :\,\supp(x)\in\supp(A)\}$.
\medskip

The next result shows that if we want to prove Conjecture~3, it is sufficient to prove it for monotone families. This will enable us to establish Theorems~8 and~7, that is, to verify the conjecture for maximal $r$-cross-intersecting pairs with a limited number of relevant
coordinates. Note that similar reduction to monotone families appears also in
\cite{FF80}.

\medskip
\noindent
{\bf Lemma 10.} {\em Let $1\le r\le n$ and let $p$ be a size vector of length
$n$.

There exists a maximal pair of $r$-cross-intersecting families $A,B\subseteq S_p$ such that both $A$ and $B$ are monotone.

If $p_i\ge3$ for all $i\in[n]$, and $A,B\subseteq S_p$ are maximal $r$-cross-intersecting families that are {\em not} balls, then there exists a pair of maximal $r$-cross-intersecting families that consists of monotone families that are not balls and have no more relevant coordinates than
$A$ or $B$.}
\medskip

\noindent
{\em Proof.} Consider the following {\em shift operations}. For any $i\in[n]$ and
$j\in[p_i]\setminus\{1\}$, for any family $A\subseteq S_p$ and any element $x\in A$, we define
\begin{align*}
\phi_i(x)&=(x_1,\dots,x_{i-1},1,x_{i+1},\dots,x_n),\\
\phi_{i,j}(x,A)&=\begin{cases}\phi_i(x)&\mbox{if }x_i=j\mbox{ and }\phi_i(x)\notin A\\x&\mbox{otherwise,}\end{cases}\\
\phi_{i,j}(A)&=\{\phi_{i,j}(x,A)\; :\; x\in A\}.
\end{align*}
Clearly, we have $|\phi_{i,j}(A)|=|A|$ for any family $A\subseteq S_p$. We claim that for
any pair of $r$-cross-intersecting families $A,B\subseteq S_p$, the families $\phi_{i,j}(A)$ and
$\phi_{i,j}(B)$ are also $r$-cross-intersecting. Indeed, if $x\in A$ and
$y\in B$ are $r$-intersecting vectors, then $\phi_{i,j}(x,A)$ and
$\phi_{i,j}(y,B)$ are also $r$-intersecting, unless $x$ and $y$ have exactly
$r$ common coordinates, one of them is $x_i=y_i=j$, and this common coordinate
gets ruined as $\phi_{i,j}(x,A)=x$ and $\phi_{i,j}(y,B)=\phi_i(y)$ (or vice
versa). However,
this is impossible, because this would imply that the vector $\phi_i(x)$ belongs to $A$, in spite of the fact that $\phi_i(x)$ and $y\in B$ are not $r$-intersecting.

If $(A,B)$ is a {\em maximal} $r$-cross-intersecting pair, then so is
$(\phi_{i,j}(A),\phi_{i,j}(B))$. When applying one of these shift operations
changes either of the families $A$ or $B$, then the total sum of all coordinates of all elements decreases. Therefore,
after shifting a finite number of times we arrive at a maximal pair of
$r$-intersecting families that cannot be changed by further shifting. We claim that
this pair $(A,B)$ is monotone. Let $y\in B$ and $y'\in S_p\setminus B$ be
arbitrary. We show that $B$ is monotone by showing that $\supp(y')$ is not
contained in $\supp(y)$. Indeed, by the maximality of the pair $(A,B)$ and using the fact that
$y'\notin B,$ there must exist $x'\in A$ such that $x'$ and $y'$ are not
$r$-cross-intersecting, and hence
$|\supp(x')\cup\supp(y')|>n-r$. Applying ``projections'' $\phi_i$ to
$x'$ in the coordinates $i\in\supp(x')\cap\supp(y)$, we obtain $x$ with
$\supp(x)=\supp(x')\setminus\supp(y)$. The shift operations $\phi_{i,j}$
do not change the family $A$, thus $A$ must be closed for the projections
$\phi_i$ and we have $x\in A$. The
supports of $x$ and $y$ are disjoint. Thus, their Hamming distance is
$|\supp(x)\cup\supp(y)|$, which is at most $n-r$, as they are
$r$-intersecting. Therefore, $\supp(x)\cup\supp(y)=\supp(x')\cup\supp(y)$ is smaller
than $\supp(x')\cup\supp(y')$, showing that $\supp(y')\not\subseteq\supp(y)$. This
proves that $B$ is monotone. By symmetry, $A$ is also monotone,
which proves the first claim of the lemma.
\smallskip

To prove the second claim, assume that $p_i\ge3$ for all $i\in[n]$. Note that
Theorem~1 establishes the lemma in the case $r=1$, so from now on we can assume without loss of generality that
$r\ge2$. Let
$A,B\subseteq S_p$ form a maximal $r$-cross-intersecting pair. By the previous
paragraph, this pair can be transformed into a monotone pair by repeated
applications of the shift
operations $\phi_{i,j}$. Clearly, these operations do not introduce new
relevant coordinates. It remains to check that the shifting operations do not
produce balls from non-balls, that is, if $A,B\subseteq S_p$ are
maximal $r$-cross-intersecting families, and $A'=\phi_{i,j}(A)$ and
$B'=\phi_{i,j}(B)$ are balls, then so are $A$ and $B$. In fact, by Theorem~4
it is sufficient to prove that one of them is a ball.
\smallskip

We saw that $A'$ and
$B'$ must also form a maximal $r$-cross-intersecting pair. Thus, by Theorem~4, there
is a set of coordinates $T\subseteq[n]$, a vector $x_0\in S_p$, and
radii $l$ and $m$ satisfying $|T|=r+l+m$ and that $A'$ and $B'$ are the
Hamming balls of radius $l$ and $m$ in coordinates $T$ around the vector
$x_0$. We can assume that $i\in T$, because otherwise $A=A'$ and we are done.
We also have that $(x_0)_i=1$, as otherwise $A'=\phi_{i,j}(A)$ is impossible.
The vectors $x\in S_p$ such that $x_i=j$ and
$$|\{k\in T\, :\, x_k\ne(x_0)_k\}|=l+1$$
are called {\em  $A$-critical}. Analogously, the vectors $y\in S_p$ such that $y_i=j$ and $$|\{k\in T\, :\,y_k\ne(x_0)_k\}|=m+1$$
are said to be {\em $B$-critical}. By the definition of
$\phi_{i,j}$, the family $A$ differs from $A'$ by including some
$A$-critical vectors $x$ and losing the corresponding vectors
$\phi_i(x)$. Symmetrically, $B\setminus B'$ consists of some $B$-critical vectors
$y$ and $B'\setminus B$ consists of the corresponding vectors $\phi_i(y)$.
Let us consider the bipartite graph $G$ whose vertices on one side are the
$A$-critical vectors $x$, the vertices on the other side are the $B$-critical
vectors $y$ (considered as disjoint families, even if $l=m$), and $x$ is
adjacent to $y$ if and only if
$|\{k\in[n]\, :\, x_k=y_k\}|=r$. If $x$ and $y$ are
adjacent, then neither the pair $(x,\phi_i(y))$, nor the pair $(\phi_i(x),y)$ is
$r$-intersecting. As $A$ and $B$ are $r$-cross-intersecting, for any pair of adjacent
vertices $x$ and $y$ of $G$, we have $x\in A$ if and only if $y\in B$.

The crucial observation is that the graph $G$ is connected. Note that this is
not the case if
$p_k=2$ for some index $k\notin T$, since all $A$-critical vectors $x$ in a
connected component of $G$ would have the same value $x_k$. However, we
assumed that $p_k>2$ for $l\in[n]$. In this case, the $A$-critical vectors $x$
and $x'$ have a common $B$-critical neighbor (and, therefore, their distance in
$G$ is $2$) if and only if the symmetric difference of the $l$ element
sets $\{k\in T\setminus\{i\}\, :\, x_k\ne(x_0)_k\}$ and $\{k\in T\setminus\{i\}\, :\,
x'_k\ne(x_0)_k\}$ have at most $2r-2$ elements. We assumed that $r>1$, so this
means that all $A$-critical vectors are indeed in the same component of the
graph $G$. Therefore, either all $A$-critical vectors belong to $A$ or none of
them does. In the
latter case, we have $A=A'$. In the former case, $A$ is the Hamming ball of
radius $l$ in coordinates $T$ around the vector $x_0'$, where $x_0'$ agrees
with $x_0$ in all coordinates but in $(x_0')_i=j$. In either case, $A$ is a
ball as required. \qed
\medskip

\noindent
{\em Proof of Theorem~8.} By Lemma~10, it is enough to restrict our attention to
monotone families $A$ and $B$. We may also assume that all coordinates are
relevant (simply drop the irrelevant coordinates). Thus, we have $n\le
r+3$.

Denote by $U_l$ the Hamming ball of radius $l$
around the all-$1$ vector in the entire set of coordinates $[n]$.
Notice that the monotone families $A$ and $B$ are $r$-cross-intersecting if and
only if for $a\in\supp(A)$ and $b\in\supp(B)$ we have $|a\cup b|\le n-r$.
We consider all possible values of $n-r$, separately.
\smallskip

If $n=r$, both families $A$ and $B$ must coincide with the singleton $U_0$.

\smallskip

If $n=r+1$, it is still true that either $A$ or $B$ is $U_0$, and hence both
families are balls. Otherwise, both
$\supp(A)$ and $\supp(B)$ have to contain at least one non-empty set, but the
union of these sets has size at most $n-r=1$, so we
have $\supp(A)=\supp(B)=\{\emptyset,\{i\}\}$ for some $i\in[n]$. This
contradicts our assumption that the coordinate $i$ is relevant for $A$.

\smallskip

If $n=r+2$, we are done if $A=B=U_1$. Otherwise, we must have a $2$-element set
either in $\supp(A)$ or in $\supp(B)$. Let us assume that a $2$-element set
$\{i,j\}$ belongs to $\supp(A)$. Then each set $b\in\supp(B)$ must satisfy
$b\subseteq\{i,j\}$. This leaves five possibilities for a non-empty monotone
family $B$, as $\supp(B)$ must be one of the following set systems:
\begin{enumerate}
\item $\{\emptyset\}$,
\item $\{\emptyset,\{i\}\}$,
\item $\{\emptyset,\{j\}\}$,
\item $\{\emptyset,\{i\},\{j\}\}$, and
\item$\{\emptyset,\{i\},\{j\},\{i,j\}\}$.
\end{enumerate}
Cases 2, 3, and 5 are not possible, because either $i$ or $j$ would not be
relevant for $B$.

In case 1, $A$ and $B$ are balls, as claimed. Nevertheless, this case is
impossible as the radii of $A$ and $B$ differ by $2$, contradicting Theorem~4.

It remains to deal with case 4. Here $\supp(A)$ consists of the sets of
size at most $1$ and the $2$-element set $\{i,j\}$. Define
$$C=\{x\in S_p\; :\; x_k=1\mbox{ for all }k\in[n]\setminus\{i,j\}\}.$$
Note that $|A|+|B|=|U_1|+|C|$, because each vector in $S_p$ appears in the same
number of sets on both sides. Thus, we have either $|A|+|B|\le2|U_1|$ or
$|A|+|B|\le2|C|$. Since $|A|>|B|$, the above inequalities imply $|A|\cdot|B|<|U_1|^2$ or $|A|\cdot|B|<|C|^2$. This
contradicts the maximality of the pair $(A,B)$, because both $U_1$ and $C$ are
$r$-intersecting. The contradiction completes the proof of the case $n-r=2$.

\smallskip

To complete the proof of Theorem~8, we need to deal with the case $n-r=3$, i.e., when
there are $r+3$ relevant coordinates. Note that, as part 1 of
Theorem~8 does not apply to this case, we have $p_i\ge 3$ for $i\in[n]$. This
slightly simplifies the following case analysis, where we consider all
containment-maximal pairs of families $(\supp(A),\supp(B))$ with the required
condition on the size of the pairwise unions.

Before considering the individual cases, we make a few simple observations.
First, we have $$\supp(A)=\{T\subseteq[n]\mid\forall U\in\supp(B):|T\cup
U\le3\},$$
$$\supp(B)=\{U\subseteq[n]\mid\forall T\in\supp(A):|T\cup U|\le3.$$
These are non-empty sets and they determine the monotone sets $A$ and
$B$.

We say that $i$ {\em dominates} $j$ in a set system $C$, if whenever $j\in T$
but $i\notin T$ for a set $T\in C$, then we have
$(T\setminus\{j\})\cup\{i\}\in C$. We say that
$i$ is {\em equivalent} to $j$ in $C$ if $i$ dominates $j$ in $C$ and $j$ also
dominates $i$. If $i$ dominates $j$ but $j$ does not dominate $i$, then
we say that $i$ {\em strictly dominates} $j$.

Note that if one of the statements {\em ``$i$ dominates $j$," ``$i$ is equivalent to $j$,"}
or {\em ``$i$ strictly dominates $j$"} holds in either $\supp(A)$ or $\supp(B)$, then
the same statement holds in both families. If $i$ strictly dominates $j$ in
$\supp(A)$, then we have $p_i\ge p_j$. Indeed, otherwise we would have
$|A'|>|A|$ and $|B'|>|B|$ for the monotone families $A'$ and $B'$ whose
supports $\supp(A')$, resp.\ $supp(B')$, are obtained from $\supp(A)$,
resp.\ $\supp(B)$, by switching the roles of $i$ and $j$. Since $A'$ and
$B'$ are $r$-cross-intersecting, this contradicts the maximality of $(A,B)$.

For equivalent coordinates $i$ and $j$ in $\supp(A)$, we may assume by symmetry
that $p_i\ge p_j$.
\smallskip

\noindent{\bf Case 1.} First assume that $\supp(A)$ contains a $3$-element
set $\{i,j,k\}$. Then all sets in $\supp(B)$ are contained in $\{i,j,k\}$.
Therefore, $\supp(B)$ must be one of the following sets, up to a suitable permutation
of the indices $i$, $j$, and $k$.
\begin{enumerate}
\item$\{\emptyset\}$,
\item$\{\emptyset,\{i\}\}$,
\item$\{\emptyset,\{i\},\{j\}\}$,
\item$\{\emptyset,\{i\},\{j\},\{k\}\}$,
\item$\{\emptyset,\{i\},\{j\},\{i,j\}\}$,
\item$\{\emptyset,\{i\},\{j\},\{k\},\{i,j\}\}$,
\item$\{\emptyset,\{i\},\{j\},\{k\},\{i,j\},\{i,k\}\}$,
\item$\{\emptyset,\{i\},\{j\},\{k\},\{i,j\},\{i,k\},\{j,k\}\}$,
\item$\{\emptyset,\{i\},\{j\},\{k\},\{i,j\},\{i,k\},\{j,k\},\{i,j,k\}\}$.
\end{enumerate}

In all of these families, $i$ dominates $j$ and $j$ dominates $k$. So we may
assume $p_i\ge p_j\ge p_k\ge3$.

Subcases 2, 5, 7, and 9 are not possible, because $i$ is not a relevant
coordinate for $A$ in them.

In subcase 1, $A$ and $B$ are balls, but, as before, this case is still
impossible, because the radii of $A$ and $B$ differ by 3.

In subcase 3, we apply Lemma~9 to the set $B$ and coordinate $i$ to obtain
$(p_i-1)^2\le p_j$, a contradiction.

In subcase 4, a similar application of Lemma~9 yields $(p_i-1)^2\le p_j+p_k-1$
with the only solution $p_i=p_j=p_k=3$. We have
$\supp(A)=\supp(U_2)\cup\{\{i,j,k\}\}$ and thus $|A|=|U_2|+8$. We further have
$|B|=7$. We must have $n\ge4$, so that $|U_1|\ge9$ and
$|U_2|\ge33$. Using these estimates, we obtain $|A|\cdot|B|<|U_1|\cdot|U_2|$, a
contradiction.

In subcase 6, we again start with Lemma~9. It yields that $(p_i-1)^2p_j\le
p_j+p_k-1$, a contradiction.

Finally, in subcase 8, we have $(p_i-1)^2(p_j+p_k-1)\le p_jp_k$ from Lemma~9, a
contradiction.
\smallskip

\noindent{\bf Case 2.} Now we assume that $\supp(A)$ contains no $3$-element
sets, but it contains two disjoint $2$-element sets $\{i,j\}$ and $\{k,l\}$. In
this case, $\supp(B)$ contains the empty set and all singletons plus one of the
following families of $2$-element sets, up to a suitable symmetry on the
indices $i$, $j$, $k$, and $l$:
\begin{enumerate}
\item$\emptyset$,
\item$\{\{i,k\}\}$,
\item$\{\{i,k\},\{i,l\}\}$,
\item$\{\{i,k\},\{j,l\}\}$,
\item$\{\{i,k\},\{i,l\},\{j,l\}\}$,
\item$\{\{i,k\},\{i,l\},\{j,k\},\{j,l\}\}$.
\end{enumerate}

Note that in subcase 1, $A$ and $B$ are balls, and subcase~6 is identical with
subcase~4 with the roles of $A$ and $B$ reversed. We
use Lemma~9 and numeric comparisons to rule out the remaining cases.

Consider the monotone ball $C$ of radius $1$ in the set of coordinates
$[n]\setminus\{i\}$. This is an $r$-intersecting family. In subcases~2 and
3, $i$ dominates all other coordinates, so we may assume that $p_i$ is maximal among
all the $p_m$ ($m\in[n]$). In both subcases considered, we have $2|C|\ge|A|+|B|$. In
subcase~2, this follows from $p_i\ge p_k$, while in subcase~3, we again need to
apply the inequality in Lemma~9. In both subcases, we have $|A|>|B|$. Thus,
$|A|\cdot|B|<|C|^2$, contradicting the maximality of the pair $(A,B)$.

A similar argument works in subcases~4 and 5. Here $i$ does not dominate $l$
(nor does it dominate $j$, in subcase~4), but it dominates all other indices, and we can
still assume by symmetry that $p_i$ is maximal. This implies that $|A|+|B|<2|C|$,
so that $|A|\cdot|B|<|C|^2$, a contradiction.
\smallskip

\noindent{\bf Case 3.} Finally, assume that Cases~1 and 2 do not hold. In this
case, for any pair $T,U\in\supp(A)$, we have $|T\cup U|\le3$ and, hence,
$\supp(B)\supseteq\supp(A)$. We can further assume by symmetry that $\supp(B)$
contains no $3$-element set and no pair of disjoint $2$-element sets. This
implies $\supp(A)=\supp(B)$ so that $A=B$. In this case, $\supp(A)$ contains the
empty set, the singletons, and a containment-maximal intersecting family of
pairs. There are only two types of such families to consider:
\begin{enumerate}
\item (star) $\supp(A)$ contains the empty set, all singletons, and all pairs
  containing some fixed coordinate $i\in[n]$.
\item (triangle) $\supp(A)$ contains the empty set, all singletons, and the pairs
  formed by two of the three distinct coordinates $i,j,k\in[n]$.
\end{enumerate}
Here subcase~1 is not possible, as $i$ is not a relevant coordinate for
$A$. In subcase~2, we may once again assume that $p_i$ is maximal. We use the
same $r$-intersecting family $C$ as in Case~2. To see that $|A|=|B|<|C|$ (a
contradiction), we use Lemma~9. \qed
\medskip

To extend Theorem~8 to somewhat larger values
of relevant coordinates (that is, to verify Conjecture~3, for instance, for the case
where there are $r+4$ relevant coordinates), we would have to go through a similar case analysis as
above. We would have to consider much more cases that correspond
to containment-maximal pairs of set systems $(U,V)$ with $|u\cup v|$ bounded
for $u\in U$ and $v\in V$. This seems to be doable, but the number of cases to
consider grows fast.
\medskip

Now we can prove our main theorem, verifying Conjecture~3 in several special cases.

\medskip
\noindent
{\em Proof of Theorem~7.} The statement about the case $p_1>r+1$ readily
follows from Lemma~2, as in this case the condition
$$\frac{2}{p_{r+1}}+\sum_{i=1}^r\frac{1}{p_i}\le1$$
holds.

To prove the other two statements in the theorem, we assume that $A$ and $B$ form a maximal $r$-cross-intersecting pair. We also assume without loss of generality that all coordinates are relevant for both families (simply drop the
irrelevant coordinates).

By Theorem~6, we have $\prod_{i=1}^rp_i\ge\prod_{i=1}^n(p_i-1)^{1-2/p_i}$, and
thus
$$\prod_{i=1}^r\frac{p_i}{(p_i-1)^{1-2/p_i}}\ge\prod_{i=r+1}^n(p_i-1)^{1-2/p_i}.$$
Here the function $x/(x-1)^{1-2/x}$ is decreasing for $x\ge3$, while $(x-1)^{1-2/x}$ is
increasing, and we have $p_i\ge p_1\ge3$. Therefore, we also have
$$\prod_{i=1}^r\frac{p_1}{(p_1-1)^{1-2/p_1}}\ge\prod_{i=r+1}^n(p_1-1)^{1-2/p_1},$$
$$p_1^r\ge(p_1-1)^{n(1-2/p_1)},$$
$$n\le\frac{r\log p_1}{(1-2/p_1)\log(p_1-1)}.$$
Simple calculation shows that the right-hand side of the last inequality is
strictly smaller than $r+4$ if $p_1\le t(r)$ for some function $t(r)=r/2+o(r)$
and, in particular, for $p_1=r+1\ge5$. In this case, we have $n\le r+3$ relevant
coordinates. Thus, Theorem~8 applies, yielding that $A$ and $B$ are balls. This
proves the last statement of Theorem~7.

\smallskip

For the proof of the second statement, note that we have already established
that $A$ and $B$ are balls in up to $r+3$ coordinates. Theorem~4 tells us that
the pair of radii must be $(0,0)$, $(0,1)$, $(1,1)$, or $(1,2)$. Simple calculation
shows that the first possibility (fixing the smallest $r$ coordinates) is
always optimal, and the cases where the two radii are unequal never yield maximal
$r$-cross-intersecting pairs. Finally, the construction with a ball of radius $1$ in $r+2$
coordinates matches the family obtained by fixing the $r$ smallest coordinates if and
only if all relevant coordinates satisfy $p_i=r+1$. This completes the proof of
Theorem~7. \qed

\section{Coordinates with $p_i=2$}\label{cr}
In many of our results, we had to assume $p_i>2$ for all coordinates of the
size vector. Here we elaborate on why the coordinates $p_i=2$ behave
differently.

For the simple characterization of the cases of equality in Theorem~1, the
assumption $k\ne 2$ is necessary. Here we characterize all maximal
cross-intersecting pairs in the case $k=2$.

Let $p=(p_1,\dots,p_n)$ be a size vector of positive integers with $k=\min_ip_i=2$ and let $I=\{i\in[n]\; :\; p_i=2\}$. For any set $W$ of functions $I\to[2]$, define the families
$$A_W=\{x\in S_p\; :\;\exists f\in W\; \mbox{such that } x_i=f(i)\; \mbox{for every } i\in I\},$$
$$B_W=\{y\in S_p\; :\;\not\exists f\in W\; \mbox{such that } y_i\ne f(i)\; \mbox{for every } i\in I\}.$$
The families $A_W$ and $B_W$ are cross-intersecting for any $W$. Moreover, if
$|W|=2^{|I|-1}$, we have  $|A_W|\cdot|B_W|=|S_p|^2/4$, so they form a maximal
cross-intersecting pair. Note
that these include more examples than just the pairs of families described in
Theorem~1, provided that $|I|>1$.

We claim that all maximal cross-intersecting pairs are of the
form constructed above. To see this, consider a maximal pair $A,B\subseteq S_p$.
We know from the proof of Theorem~1 that $x\in A$ if and only if $x'\notin B$, where $x'$
is defined by $x'_i=(x_i+1\bmod p_i)$ for all $i\in[n]$.  Let $j\in[n]$ be a
coordinate with $p_j>2$. By the same argument, we also have that
$x\in A$ holds if and only if $x''\notin B$, where $x''_i=x'_i$ for
$i\in[n]\setminus\{j\}$ and $x''_j=(x_j+2\bmod p_j)$. Thus, both $x'$ and
$x''$ belong to $B$ or neither of them does. This holds for every vector $x'$,
implying that $j$ is irrelevant for the family $B$ and thus also for $A$.

As there are no relevant coordinates for $A$ and $B$ outside the set $I$ of
coordinates with $p_i=2$, we can choose a set $W$ of functions from $I$ to
$[2]$ such that $A=A_W$. This makes
$$B=\{y\in S_p\;:\;y\mbox{ intersects all } x\in A\}=B_W.$$
We have $|A|+|B|=|S_p|$ and $|A|\cdot|B|=|S_p|^2/4$ if and only if
$|W|=2^{|I|-1}$.
\smallskip

The size vector $p=(2,\dots,2)$ of length $n$ is well studied. In this case, $S_p$ is the $n$-dimensional hypercube. If $r>1$, then all maximal
$r$-cross-intersecting pairs have an unbounded number of relevant coordinates,
as a function of $n$. Indeed, the density $|A|\cdot|B|/|S_p|^2$ is at most
$1/4$ for cross-intersecting pairs $A,B\subseteq S_p$, and strictly less than
$1/4$ for $r$-cross-intersecting families if $r>1$. Furthermore, if the number of
relevant coordinates is bounded, then this density is bounded away from $1/4$,
while if $A=B$ is the ball of radius $(n-r)/2$ in all the coordinates, then
the same density approaches $1/4$.

One can also find many maximal $2$-cross-intersecting pairs that are not
balls. For example, in the $3$-dimensional hypercube the families
$A=\{0,0,0),(0,1,1)\}$ and $B=\{(0,0,1),(0,1,0)\}$ form a maximal
$2$-cross-intersecting pair.

\smallskip
Finally, we mention that there is a simple connection between the problem
discussed in this paper and a question related to communication
complexity. Consider the following two-person communication game: Alice and
Bob each receive a vector from $S_p$, and they have to decide whether the
vectors are $r$-intersecting. In the {\em communication matrix} of such a game, the
rows are indexed by the possible inputs of Alice, the columns by the possible
inputs of Bob, and an entry of the matrix is $1$ or $0$ corresponding to the
``yes'' or ``no'' output the players have to compute for the corresponding
inputs. In the study of communication games, the submatrices of this
matrix in which all entries are equal play a special role. The largest
area of an all-$1$ submatrix is
the maximal value of $|A|\cdot|B|$ for $r$-cross-intersecting families
$A,B\subseteq S_p$.

\bigskip
\noindent{\bf Acknowledgment.} We are indebted to G. O. H. Katona, R. Radoi\v ci\'c, and D. Scheder for their valuable remarks, and to an anonymous referee for calling our attention to the manuscript of Borg~\cite{Bo14}.

%\begin{acknowledgements}
%\end{acknowledgements}


\begin{thebibliography}{}



%\small
\bibitem{AhK98}
R. Ahlswede and L. H. Khachatrian, The diametric theorem in Hamming spaces--optimal anticodes, {\em Adv. in Appl. Math.}, {\bf 20}, 429--449 (1998).

\bibitem{Be74}
C. Berge, Nombres de coloration de l'hypegraphe $h$-parti complet, {\em Hypergraph Seminar, Lecture Notes in Math.}, Springer-Verlag, Heidelberg, {\bf 411}, 13--20 (1974).

%\bibitem[Be05]{Be05} C. Bey: On cross-intersecting families of sets, {\em Graphs Combin.} {\bf 21} 2005), 161--168.

\bibitem{Bol86}
B. Bollob\'as, {\em Combinatorics. Set Systems, Hypergraphs, Families of Vectors and Combinatorial Probability}, Cambridge University Press, Cambridge (1986).

%\bibitem[BoL97]{BoL97} B. Bollob\'as and I. Leader: An Erd\H os-Ko-Rado theorem for signed sets, {\em Comput. Math. Applic.} {\bf 34} (1997), 9--13.

\bibitem{Bo08}
P. Borg, Intersecting and cross-intersecting families of labeled sets, {\em Electron. J. Combin.}, {\bf 15} N9 (2008).

%\bibitem[Bo09a]{Bo09a} P. Borg: On $t$-intersecting families of signed sets and permutations, {\em Discrete Math.} {\bf 309} (2009), 3310--3317.

%\bibitem[Bo09b]{Bo09b} P. Borg: Extremal $t$-intersecting sub-families of hereditary families, {\em J. London Math. Soc.} {\bf 79} (2009), 167--185.

%\bibitem[Bo09c]{Bo09c} P. Borg: A short proof of a cross-intersection theorem of Hilton, {\em Discrete Math.} {\bf 309} (2009), 4750--4753.

%\bibitem[Bo10]{Bo10} P. Borg: Cross-intersecting families of permutations, {\em J. Combin. Theory Ser. A} {\bf 117} (2010), 483--487.

\bibitem{Bo14}
P. Borg, Cross-intersecting integer sequences, preprint, arXiv:1212.6955.

%\bibitem[BoL10]{BoL10} P. Borg and I. Leader: Multiple cross-intersecting families of signed sets, {\em J. Combin. Theory Ser. A} {\bf 117} (2010), 583--588.

%\bibitem[CaK03]{CaK03} P. J. Cameron and C. Y. Ku: Intersecting families of permutations, {\em European J. Combin.} {\bf 24} (2003), 881--890.

%\bibitem[DeF77]{DeF77} M. Deza and P. Frankl: On the maximum number of permutations with given maximal or minimal distance, {\em J. Combin. Theory Ser. A} {\bf 22} (1977), 352--360.

\bibitem{DeF83}
M. Deza and P. Frankl, Erd\H os-Ko-Rado theorem--22 years later, {\em SIAM J. Alg. Disc. Methods}, {\bf 4}, 419--431 (1983).

%\bibitem[EFP11]{EFP11} D. Ellis, E. Friedgut, and H. Pilpel: Intersecting families of permutations, {\em J. Amer. Math. Soc.} {\bf 24} (2011), 649--682.

\bibitem{EKR61}
P. Erd\H os, C. Ko, and R. Rado, Intersection theorems for systems of finite sets, {\em Quart. J. Math. Oxford Ser. 2}, {\bf 12}, 313--318 (1961).

\bibitem{F95}
P. Frankl, Extremal set systems, in {\em Handbook of Combinatorics} (R. Graham et al., eds.), Elsevier, Amsterdam, 1293--1329 (1995).

\bibitem{FF80}
P. Frankl and Z. F\"uredi, The Erd\H os-Ko-Rado theorem for integer sequences, {\em SIAM J. Algebraic Discrete Methods}, {\bf 1}, 376--381 (1980).

\bibitem{FLST14}
P. Frankl, S. J. Lee, M. Siggers, and N. Tokushige, An Erd\H{o}s-Ko-Rado theorem for cross-intersecting families, {\em J. Combin. Theory Ser. A}, {\bf 128}, 207--249 (2014).

\bibitem{FT99}
P. Frankl and N. Tokushige, The Erd\H os-Ko-Rado theorem for integer sequences, {\em Combinatorica}, {\bf 19}, 55-63 (1999).

%\bibitem[FW86]{FW86} P. Frankl and R. M. Wilson: The Erd\H os-Ko-Rado theorem for vector spaces, {\em J. Combin. Theory Ser. A} {\bf 43} (1986), 228--236.

%\bibitem[GeWZ12]{GeWZ12} X. B. Geng, J. Wang, and H. J. Zhang: Structures of independent sets in direct products of some vertex-transitive graphs, {\em Acta Math. Sin. (Engl. Ser.)} {\bf 28} (2012): 697--706.

\bibitem{GK78}
C. Greene and D. J. Kleitman, Proof techniques in the ordered sets, in: {\em Studies in Combinatorics}, Math. Assn. America, Washington DC, 22--79 (1978).

\bibitem{Hi77}
A. J. W. Hilton, An intersection theorem for a collection of families of subsets of a finite set, {\em J. London Math. Soc.}, {\bf 2}, 369--384 (1977).

%\bibitem[H75]{H75} W. N. Hsieh: Intersection theorems for systems of finite vector spaces, {\em Discrete Math.} {\bf 12} (1975), 1--16.

\bibitem{K72}
G. O. H. Katona, A simple proof of the Erd\H os-Ko-Rado theorem, {\em J. Combin. Theory Ser. B}, {\bf 13}, 183--184 (1972).

%\bibitem[LW07]{LW07} Y. S. Li and J. Wang: Erd\H os-Ko-Rado-type theorems for colored sets, {\em Electron. J. Combin.} {\bf 14(1)} (2007), N1.

\bibitem{Liv79}
M. L. Livingston, An ordered version of the Erd\H os-Ko-Rado theorem, {\em J. Combin. Theory Ser. A}, {\bf 26}, 162--165 (1979).

%\bibitem[MaT89] M. Matsumoto and N. Tokushige: The exact bound in the Erd\H os-Rado-Ko theorem for cross-intersecting families, {\em J. Combin. Theory Ser. A} {\bf 52} (1989), 90--97.

\bibitem{Mo82}
A. Moon, An analogue of the Erd\H os-Ko-Rado theorem for the Hamming schemes H(n,q), {\em J. Combin. Theory Ser. A}, {\bf 32}, 386--390 (1982).

\bibitem{Py86}
L. Pyber, A new generalization of the Erd\H os-Rado-Ko theorem, {\em J. Combin. Theory Ser. A}, {\bf 43}, 85--90 (1986).

\bibitem{St80}
D.  Stanton, Some Erd\H os-Ko-Rado theorems for Chevalley groups, {\em SIAM J. Algebraic Discrete Methods} {\bf 1}, 160--163 (1980).

%\bibitem[To10]{To10} N. Tokushige, On cross $t$-intersecting families of sets, {\em J. Combin. Theory Ser. A} {\bf 117} (2010), 1167--1177.

\bibitem{To13}
N. Tokushige, Cross $t$-intersecting integer sequences from weighted Erd\H os-Ko-Rado, {\em Combin. Prob. Comput.}, {\bf 22}, 622--637 (2013).

%\bibitem[WaZ08]{WaZ08} J. Wang and S. J. Zhang: An Erd\H os-Ko-Rado-type theorem in Coxeter groups, {\em European J. Combin.} {\bf 29} (2008), 1112--1115.

%\bibitem[WaZ11a]{WaZ11a} J. Wang and H. J. Zhang: Intersecting families in a subset of boolean lattices, {\em Electron. J. Combin.} {\bf 18} (2011), P237.

%\bibitem[WaZ11b]{WaZ11b} J. Wang and H. J. Zhang: Cross-intersecting families and primitivity of symmetric systems, {\em J. Combin. Theory Ser. A} {\bf 118} (2011), 455--462.

%\bibitem[WaZ13]{WaZ13} J. Wang and H. J. Zhang: Nontrivial independent sets of bipartite graphs and cross-intersecting families, {\em J. Combin. Theory Ser. A} {\bf 120} (2013), 129--141.

\bibitem{Zh13}
H. Zhang, Cross-intersecting families of labeled sets, {\em Electron. J. Combin.} {\bf 20(1)}, P17 (2013).




\end{thebibliography}


\end{document}